\documentclass{amsart}

\usepackage{amsmath,amssymb}
\usepackage{dsfont,amsthm}
\usepackage{graphicx}
\newcommand{\N}{{\mathds N}}
\newcommand{\R}{{\mathds R}}

\newcommand{\Var}{{\rm Var}}

\newtheorem{lemma}{Lemma}[section]

\newtheorem{theorem}{Theorem}[section]

\begin{document}
\title{Bootstrap for $U$-Statistics: A new approach}
   \author{Olimjon Sh. Sharipov, Johannes Tewes, Martin Wendler}
   
  \date{\today}

\begin{abstract}
Bootstrap for nonlinear statistics like $U$-statistics of dependent data has been studied by several authors. This is typically done by producing a bootstrap version of the sample and plugging it into the statistic. We suggest an alternative approach of getting a bootstrap version of $U$-statistics, which can be described as a compromise between bootstrap and subsampling. We will show the consistency of the new method and compare its finite sample properties in a simulation study.
\end{abstract}

\maketitle

\section{Introduction}
In many statistical applications, the asymptotic limit distribution cannot be used to construct tests or confidence intervals, as the limit might be unknown or dependent on unknown parameters, which are hard to estimate. Bootstrap versions of the statistical estimator provide a nonparametric alternative, so establishing the consistency of such bootstrap procedure is important. In the case of independent and identically distributed (iid) random variables, the main idea of the bootstrap consists of replacing the original sample $X_{1},...,X_{n}$ of observations with unknown marginal distribution function $F(x)$ by a new iid sample $X_{1}^{*},...,X_{n}^{*}$ with the marginal distribution function $F_{n}(x)$ which is the empirical distribution function constructed from the original sample (see Efron \cite{efro}). As a bootstrap version of a statistic $T(X_{1},...,X_{n})$ one can take the plug-in version $T(X_{1}^{*},...,X_{n}^{*})$.
The validity of bootstrap for the sample mean of iid observations was first established by Bickel and Freedman \cite{bick} and Singh \cite{sing}. It is well-known that in some cases the bootstrap method provides a better approximation to the distribution of the statistic than, for instance, normal approximation (see Hall \cite{hall}). Especially when sample size is relatively small, bootstrap has some preferences.

In the case of dependent observations this idea of getting a new sample does not work since a new iid sample does not capture the dependence structure, see Singh \cite{sing}. Therefore several so-called blockwise bootstrap methods of getting a new sample were introduced (for overview of blockwise bootstrap methods see Lahiri \cite{lahi}). In all these methods, we resample from blocks of $l$ consecutive observations, so that inside the blocks the dependence structure will be kept. We will consider two of them: The circular block bootstrap method introduced by Politis and Romano \cite{poli} and nonoverlapping block bootstrap method by Carlstein \cite{carl}. In the case of circular block bootstrap method we extend the original sample to $ X_{1},...,X_{n},X_{1},...,X_{l-1}$ where $l=l(n)$ is a block length. To get a new sample $X_{1,1}^{*},...,X_{1,ml}^{*}$ we choose randomly and independently $l$ consecutive observations of the sample $m=[\frac{n}{l}]$ times with
\begin{multline*}
P^{\star}(X_{1,kl+1}^{*}=X_{j}, ... , X_{1,(k+1)l}^{*}=X_{j+l-1}):=\\
P\left(X_{1,kl+1}^{*}=X_{j}, ... , X_{1,(k+1)l}^{*}=X_{j+l-1}\big|X_1,\ldots,X_n\right)=\frac{1}{n} 
\end{multline*}
for $j=1,...,n ,k=0,...,m-1.$ Here and in what follows we denote by $P^{\star}, E^{\star}, \Var^{\star}$ the conditional probability, conditional expectation and conditional variance respectively. For instance the bootstrap version of sample mean $\overline{X}_{n}$ will be
\begin{equation*}
 \overline{X}_{n}^{1*}=\dfrac{1}{ml}\sum_{i=1}^{ml}X_{1,i}^{*}.
\end{equation*}
In the case of nonoverlapping block bootstrap method we construct a new sample $X_{2,1}^{*},...,X_{2,ml}^{*}$ by choosing randomly and independently blocks of $l$ consecutive observations of the sample $ X_{1},...,X_{n}$ as above $m=[\frac{n}{l}]$ times with
\begin{equation*}
P^{\star}(X_{2,kl+1}^{*}=X_{(j-1)l+1}, ... , X_{2,(k+1)l}^{*}=X_{jl})=\frac{1}{m} 
\end{equation*}
for $j=1,...,m ,k=0,...,m-1$. In this case the bootstrap version of the sample mean is
\begin{equation*}
 \overline{X}_{n}^{2*}=\dfrac{1}{ml}\sum_{i=1}^{ml}X_{2,i}^{*}.
\end{equation*}
For weakly dependent observations strong consistency (almost surely convergence) of the $\overline{X}_{n}^{1*}$ and $\overline{X}_{n}^{2*}$ were proved in Shao and Yu \cite{shao} and Peligrad \cite{peli}. Gon\c{c}alves and White \cite{gonc} and Dehling, Sharipov and Wendler \cite{deh3} extended these results to functionals of mixing processes.
In this paper we will concentrate on bootstrap for $U$-statistics of weakly dependent observations. Bootstrap for $U$-statistics of independent observations were studied by Bickel and Freedman \cite{bick}, Arcones and Gine \cite{arco}, Dehling and Mikosch \cite{deh1}, Leucht and Neumann \cite{leuc}. Recently, Dehling and Wendler \cite{dehl}, Sharipov and Wendler \cite{shar} and Leucht and Neumann \cite{leu2} established the consistency of bootstrap estimators for $U$-statistics of weakly dependent observations. In the aforementioned papers the same idea of getting bootstrap versions of $U$-statistics have been used, which is based on the plug-in principle.
We will now introduce $U$-statistics. Let $(X_{n})_{n\in\N}$ be a stationary sequence of random variables with a common distribution function $F(x)$. For simplicity reasons we consider $U$-statistics of degree two i.e.:
\begin{equation*}
 U_{n}(h)=\frac{2}{n(n-1)}\sum_{1\leq i<j\leq n} h(X_{i},X_{j}) 
\end{equation*}
where $h:R\times R\rightarrow R $ is the symmetric and measurable kernel. The kernel $h(x,y)$ is called degenerate, if 
\begin{equation*}
 Eh(X_{1},x)=0 \ \  \text{for all} \ x\in R . 
\end{equation*}
If this does not hold we call the kernel $h(x,y)$ and the corresponding $U$-statistics nondegenerate. In what follows we will only consider $U$-statistics with a nondegenerate kernel. In the case of iid observations $ X_{1},...,X_{n}$ Bickel and Freedman \cite{bick} used the bootstrap sample $X_{1}^{*},...,X_{n}^{*} $ of conditional independent observations with common distribution function $F_{n}(x)$, which is the empirical distribution function constructed by original sample. In order to get bootstrap versions of $U$-statistics they plugged the bootstrapped observations in $U$-statistics i.e.
\begin{equation*}
U_{n}^{*}(h)=\frac{2}{n(n-1)}\sum_{1\leq i<j\leq n} h(X_{i}^{*},X_{j}^{*})   
\end{equation*}
In  the case of dependent observations, the same idea was explored by Dehling, Wendler \cite{dehl} using the circular block bootstrap method, the corresponding bootstrap version is
\begin{equation*}
U_{n,1}^{*}(h)=\frac{2}{ml(ml-1)}\sum_{1\leq i<j\leq ml} h(X_{1,i}^{*},X_{1,j}^{*}),
\end{equation*}
while in the case of nonoverlapping block bootstrap method the corresponding bootstrap version of $U$-statistics is
\begin{equation*}
U_{n,2}^{*}(h)=\frac{2}{ml(ml-1)}\sum_{1\leq i<j\leq ml} h(X_{2,i}^{*},X_{2,j}^{*}).
\end{equation*}
Consistency of $U_{n}^{*}(h)$ (and $U_{n,k}^{*}(h), k=1,2$ ) can be proved using the Hoeffding decomposition (see \cite{hoef}):
\begin{equation*}
U_{n}^{*}(h)=\theta +\frac{2}{n}\sum_{i=1}^{n} h_{1}(X_{i}^{*})+\frac{2}{n(n-1)}\sum_{1\leq i<j\leq n} h_{2}(X_{i}^{*},X_{j}^{*})   
\end{equation*}
where
\begin{align*}
\theta&:=Eh\left(X_1,X_2\right),\\
h_1(x)&:=Eh(x,X_{2})-\theta, \\
h_2(x,y)&:=h(x,y) - h_1(x) -h_1(y) -\theta.
\end{align*}
In order to prove consistency of $U_{n}^{*}(h)$ it is enough to show the convergence of bootstrapped distribution of the sample mean for the second summand of this decomposition and convergence to zero (in probability or almost surely) of the third summand. The goal of this paper is to suggest a new resampling method for $U$-statistics. The main idea is the following: We suggest to draw with replacement from the $U$-statistics calculated on subsamples instead of drawing with replacement from blocks of observations and then plug them in $U$-statistics. We introduce subsets $\{1,\ldots,n\}^2$ by
\begin{align*}
B_{1}^{'}(i)=\lbrace (a,b): i\leq a,b \leq i+l-1, a<b \rbrace \ \ i=1,...,n ,\\
B_{2}^{'}(i)=\lbrace (a,b): l(i-1)+1 \leq a,b \leq il, a<b \rbrace \ \ i=1,...,m=[\frac{n}{l}]
\end{align*}
where $l=l(n)\rightarrow\infty$ us called block length. Furthermore, let $\lbrace T_{1}(j),j\geq 1\rbrace$ and $\lbrace T_{2}(j),j\geq 1\rbrace$ be two sequences of iid random variables with distributions:
\begin{align*}
P(T_{1}(1)=i)=\frac{1}{n}, \ \ i=1,...,n ,\\
P(T_{2}(1)=i)=\frac{1}{m}, \ \ i=1,...,m.
\end{align*}
A fixed realisations of the sample $X_{1},...,X_{n}$ we denote by $x_{1},...,x_{n}$. We set
\begin{align*}
\overline{x}_{1}(i)&=\frac{2}{l(l-1)}\sum_{(a,b)\in B^{'}_{1}(i)} h(x_{a},x_{b}),\ \ \ i=1,...,n , \\ 
\overline{x}_{2}(i)&=\frac{2}{l(l-1)}\sum_{(a,b)\in B^{'}_{2}(i)} h(x_{a},x_{b}),\ \ \ i=1,...,m , \\
u_{j}^{1*}&=\sum_{i=1}^{n}I(T_{1}(j)=i)\overline{x}_{1}(i)=\overline{x}_{1}(T_{1}(j)),\\
u_{j}^{2*}&=\sum_{i=1}^{m}I(T_{2}(j)=i)\overline{x}_{2}(i)=\overline{x}_{2}(T_{2}(j)),
\end{align*}
where $I(\cdot)$ is the indicator function. So $\overline{x}_{1}(i)$ respectively $\overline{x}_{2}(i)$ are the $U$-statistics calculated from the $i$-th block and $u_{j}^{1*}$ and $u_{j}^{2*}$ are the results of the drawing with replacement from these $U$-statistics. In the case of $\overline{x}_{1}(i)$ we assume that $x_{k}=x_{k-n}$ if $k>n$. As the bootstrap versions of $U$-statistics $U_n(h)$ we take the following:
\begin{align*}
U_{n}^{1*}(h)&=\frac{1}{m}\sum_{i=1}^{m}u_{j}^{1*} ,\\
U_{n}^{2*}(h)&=\frac{1}{m}\sum_{i=1}^{m}u_{j}^{2*}.
\end{align*}
Note that
\begin{align*}
E^{\star}U_{n}^{1*}(h)=\frac{1}{nl(l-1)}\sum_{i=1}^{n}\sum_{(a,b)\in B^{'}_{1}(i)} h(x_{a},x_{b})\\
E^{\star}U_{n}^{2*}(h)=\frac{1}{ml(l-1)}\sum_{i=1}^{m}\sum_{(a,b)\in B^{'}_{2}(i)} h(x_{a},x_{b}).
\end{align*}
In the next sections we will state our main result: the strong consistency of the bootstrap based on $U_{n}^{k*}(h)$, $k=1,2$. The new approach reduces the computational burden of the Monte Carlo method usually used for the bootstrap. After calculating the values $u_{j}^{k*}$, for every run of the Monte Carlo evaluation one need only $O(\frac{n}{l})$ calculation steps, while for the plug-in method, in every run the $U$-statistic has to be calculated again in $O(n^2)$ steps. As the linear parts of the plug-in bootstrap version and the new bootstrap version of a $U$-statistic are the same, we expect a similar behaviour.

We will give an upper bound for the mean square error (MSE) of the variances of $U_{n}^{k*}(h)$, $k=1,2$, which suggest a choice of the block length of order $O(\sqrt{n})$. In section 3, we will present our simulation results. We consider the sample variance as a $U$-statistic and we will compare the new bootstrap approach with plug-in bootstrap and subsampling. The proofs of the main results will be given in section 4.
 
\section{Main results}
Let $\left(X_n\right)_{n\in\mathds{N}}$ be a stationary sequence of random variables with values in a separable linear space. We will assume that this sequence satisfies some form of short range dependence. Namely we will consider strong mixing and absolute regularity conditions. Recall that strong mixing coefficients $(\alpha(k))_{k\in\N}$ and absolute regularity coefficients $(\beta(k))_{k\in\N}$ are defined as
\begin{align*}
\alpha (n)=\sup\lbrace\vert P(AB)-P(A)P(B)\vert : A \in {\mathcal {F}}_{1}^{k},B\in {\mathcal {F}}_{k+n}^{\infty},k\in \N \rbrace ,\\
\beta (n)=\sup\lbrace E(\sup\vert P(B/{\mathcal {F}}_{1}^{k})-P(B)\vert : B\in {\mathcal {F}}_{k+n}^{\infty},k\in \N \rbrace 
\end{align*}
where ${\mathcal {F}}_{a}^{b}$ is the $\sigma$ - field generated by $ X_{a},...,X_{b}$. Now we can formulate our results.

\begin{theorem}\label{theo1} Let $\left(X_n\right)_{n\in\mathds{N}}$ be a stationary sequence of absolutely regular random variables. Assume that the following conditions hold (for some $\delta >0$)
\begin{itemize}
\item $\iint\left|h\left(x_{1},x_{2}\right)\right|^{2+\delta}dF\left(x_{1}\right)dF\left(x_{2}\right)\leq M$ for some $M<\infty$,
\item $E\left|h\left(X_{1},X_{1+k}\right)\right|^{2+\delta} \leq M$ for all $k\in\N$,
\item $\beta\left(n\right)=O\left(n^{-\rho}\right)$ for some $\rho>\frac{2+\delta}{\delta}$,
\item $l (n) \leq C n^\epsilon$ for some $\epsilon\in(0,1)$.
\end{itemize}
Then as $n\rightarrow\infty$ for $i=1,2$
\begin{equation*}
 \sup_{x\in\mathds{R}}\left|P^{\star}\left[\sqrt{ml}\left(U^{i\star}_{n}\left(h\right)-E^{\star}\left[U^{i\star}_{n}\right]\right)\leq x\right]-P\left[\sqrt{n}\left(U_{n}\left(h\right)-\theta\right)\leq x\right]\right|\rightarrow0.
\end{equation*}
\begin{equation*}
\left|\Var^\star\left[\sqrt{ml}U_n^{i\star}\left(h\right)\right]-\Var\left[\sqrt{n}U_n\left(h\right)\right]\right|\rightarrow0,
\end{equation*}
in probability.
\end{theorem}
In the next theorem we consider strongly mixing random variables. This is a weaker assumption on the dependence, but in this case we assume that the kernel satisfies the following condition: We say that a kernel is $\mathcal{P}$-Lipschitz-continuous, if there exists a constant $C>0$ such that 
\begin{equation*}
 E\left[\left|h\left(X,Y\right)-h\left(X',Y\right)\right|I\left(\left|X-X'\right|\leq\epsilon\right)\right]\leq C\epsilon
\end{equation*}
for every $\epsilon>0$, every pair $X$ and $Y$ with the common distribution $\mathcal{P}_{X_{1},X_{k}}$ for some $k\in\mathds{N}$ or $\mathcal{P}_{X_{1}}\times\mathcal{P}_{X_{1}}$ and $X'$ and $Y$ also with one of these common distributions. For the examples of kernels which satisfy above condition see \cite{dehl}, \cite{shar}.
\begin{theorem}\label{theo2}
Let $\left(X_n\right)_{n\in\mathds{N}}$ be a stationary sequence of strongly mixing random variables. Assume that the following conditions hold
\begin{itemize}
\item $\iint\left|h\left(x_{1},x_{2}\right)\right|^{2+\delta}dF\left(x_{1}\right)dF\left(x_{2}\right)\leq M$,
\item $E\left|h\left(X_{1},X_{1+k}\right)\right|^{2+\delta} \leq M$ for all $k\in\N$,
\item $l (n) \leq C n^\epsilon$ for some $\epsilon\in(0,1)$,
\item $h$ is $\mathcal{P}$-Lipschitz-continuous,
\item $E\left|X_{1}\right|^{\gamma}>0$ for some $\gamma>0$,
\item $\alpha(n)=O(n^{-\rho})$ for some $ \rho>\frac{3\gamma\delta+\delta+5\gamma+2}{2\gamma\delta}$.
\end{itemize}
Then the statements of Theorem \ref{theo1} hold.
\end{theorem}
In the next theorem we give bounds for the mean squared error to give a deeper insight into the properties of bootstrap versions of $U$-statistics.
\begin{theorem}\label{theo3} Let $\left(X_n\right)_{n\in\mathds{N}}$ be a stationary sequence of absolutely regular random variables. Assume that the following conditions hold:
\begin{itemize}
\item $\iint\left|h\left(x_{1},x_{2}\right)\right|^{12+\delta}dF\left(x_{1}\right)dF\left(x_{2}\right)\leq M$ for some $\delta >0$,
\item $E\left|h\left(X_{1},X_{1+m}\right)\right|^{4+\frac{2}{3}\delta} \leq M$ for all $m\in\N$,
\item $\beta\left(n\right)=O\left(n^{-\frac{3(6+\delta^{'})}{\delta^{'}}}\right)$ for some $0<\delta^{'}<\delta$.
\end{itemize}
Then for $k=1,2$
\begin{equation*}
E\left(\Var^{\star}(\sqrt{ml}U_{n}^{k*}(h))-\Var(\sqrt{n}U_{n}(h))\right)^{2}=O(\frac{1}{l}+\frac{l}{n}).
\end{equation*}
\end{theorem}
Choosing a block length of order $1/\sqrt{n}$, we can achieve that the mean squared error of $\Var^{\star}(\sqrt{ml}U_{n}^{k*}(h))$ is of order $O(1/\sqrt{n})$.

\section{Simulation}
We consider the estimator $\sigma_n^2$ for the variance $\sigma^{2}=Var(X_{1})$ which is a $U$-statistic with the kernel $h(x,y)=\dfrac{1}{2}(x-y)^{2}$:
\begin{equation*}
\sigma_{n}^{2}=\dfrac{1}{n-1}\sum^n_{i=1}(X_{i}-\overline{X}_{n})^{2}=\frac{2}{n(n-1)}\sum_{1\leq i<j \leq n}\frac{1}{2}(X_i-X_j)^2.
\end{equation*}
$(X_n)_{n\in\N}$ is a stationary, Gaussian, autoregressive process with $X_n=\alpha X_{n-1}+\epsilon_n$, where $(\varepsilon_{n})_{n\in\N}$  is a sequence of iid standard normal random variables and $\alpha\in\{0.2,0.4,0.6\}$. 
We will compare three methods for constructing confidence intervals:
\begin{itemize}
\item The circular plug-in bootstrap (already used in Dehling, Wendler \cite{dehl}).
\item The new bootstrap version $U_n^{1\star}$, the drawing with replacement from $U$-statistics caculated on subsamples.
\item Subsampling: The estimator of the unkown distribution is given by the empirical distribution function of the $U$-statistics caculated on subsamples, see Politis and Romano \cite{pol3}.
\end{itemize}
For each combination of construction method, AR-coefficient, sample size ($n=50,100,$ $200,400$) and block length, we have simulated 10.000 samples and evaluated the empirical probability of the 95\% confidence interval to cover the true parameter. The two bootstrap version are evaluated by the Monte Carlo method with 1000 times drawing with replacement. The results are summarized in the table below.

\begin{table}
\begin{center}
\begin{tabular}{lc||c|c|c|}
\hline
&$\alpha$ & 0.2 & 0.4 & 0.6\\
$n$&$l$  & &&\\ 
\hline \hline 
50 &3& 0.876/0.895/0.827 &  0.863/0.791/0.693 & 0.794/0.600/0.569 \\ 
 &5& 0.865/0.884/0.823 & 0.838/0.819/0.737 & 0.765/0.677/0.583\\ 
• & 7& 0.862/0.871/0.814 & 0.823/0.818/0.753 & 0.762/0.718/0.614\\ 
• & •10&  0.835/0.849/0.787 & 0.808/0.814/0.743 & 0.757/0.730/0.607\\ 
\hline 
100 &3& 0.908/0.922/0.823 & 0.894/0.817/0.730 & 0.852/0.600/0.507\\ 
 &5& 0.911/0.910/0.852 & 0.896/0.838/0.770 & 0.851/0.682/0.617\\ 
• &7& 0.909/0.909/0.849 & 0.874/0.857/0.793 & 0.853/0.737/0.676\\ 
• &10& 0.890/0.902/0.851 & 0.880/0.859/0.815 & 0.849/0.767/0.716\\ 
\hline 
200 &5&0.926/0.930/0.873& 0.908/0.859/0.813& 0.883/0.684/0.620 \\ 
 &7& 0.925/0.926/0.881& 0.905/0.870/0.831& 0.891/0.738/0.704 \\ 
• &10& 0.924/0.920/0.889& 0.911/0.885/0.857& 0.887/0.805/0.761\\ 
• &15& 0.910/0.910/0.884& 0.902/0.893/0.849& 0.885/0.835/0.800 \\ 
\hline 
400 &7& 0.927/0.932/0.896 & 0.927/0.887/0.855 & 0.920/0.748/0.711\\ 
 &10& 0.931/0.927/0.906 & 0.924/0.901/0.873 & 0.920/0.814/0.777\\ 
 &15& 0.923/0.929/0.901 & 0.923/0.904/0.882 & 0.902/0.847/0.828\\ 
•&20& 0.924/0.918/0.906 & 0.902/0.909/0.887 & 0.904/0.866/0.849\\ 
\hline 
\end{tabular} 
\caption{Simulated coverage propabilites for plug-in bootstrap/ new bootstrap/ subsampling, sample size $n=50,100,200,400$, AR-coefficient $\alpha=0.2,0.4,0.6$, block length $l=3,5,7,10,15,20$}
\end{center}
\end{table}
In general, we draw the following conclusions from the simulation study: All three methods lead to confidence intervals with a coverage probabilty lower than the nominal confidence level of 95\%. The performance of the subsampling is worse than the performance of two bootstrap methods in all situations. If the dependence in the AR-process is weak ($\alpha=0.2$), both bootstrap methods lead to comparable results, while for stronger dependence ($\alpha=0.4,0.6$), the plug-in bootstrap has a better performance.

For the plug-in bootstrap, the block length can be choosen smaller than for the new bootstrap method or for subsampling, especially in the case of stronger dependence ($\alpha=0.4,0.6$). In Figur 1 below, we give a more detailed picture of the coverage probabilities for different block lengths in the case $n=200$, $\alpha=0.4$.
\begin{figure}
        \centering
        \includegraphics[width=0.8\textwidth]{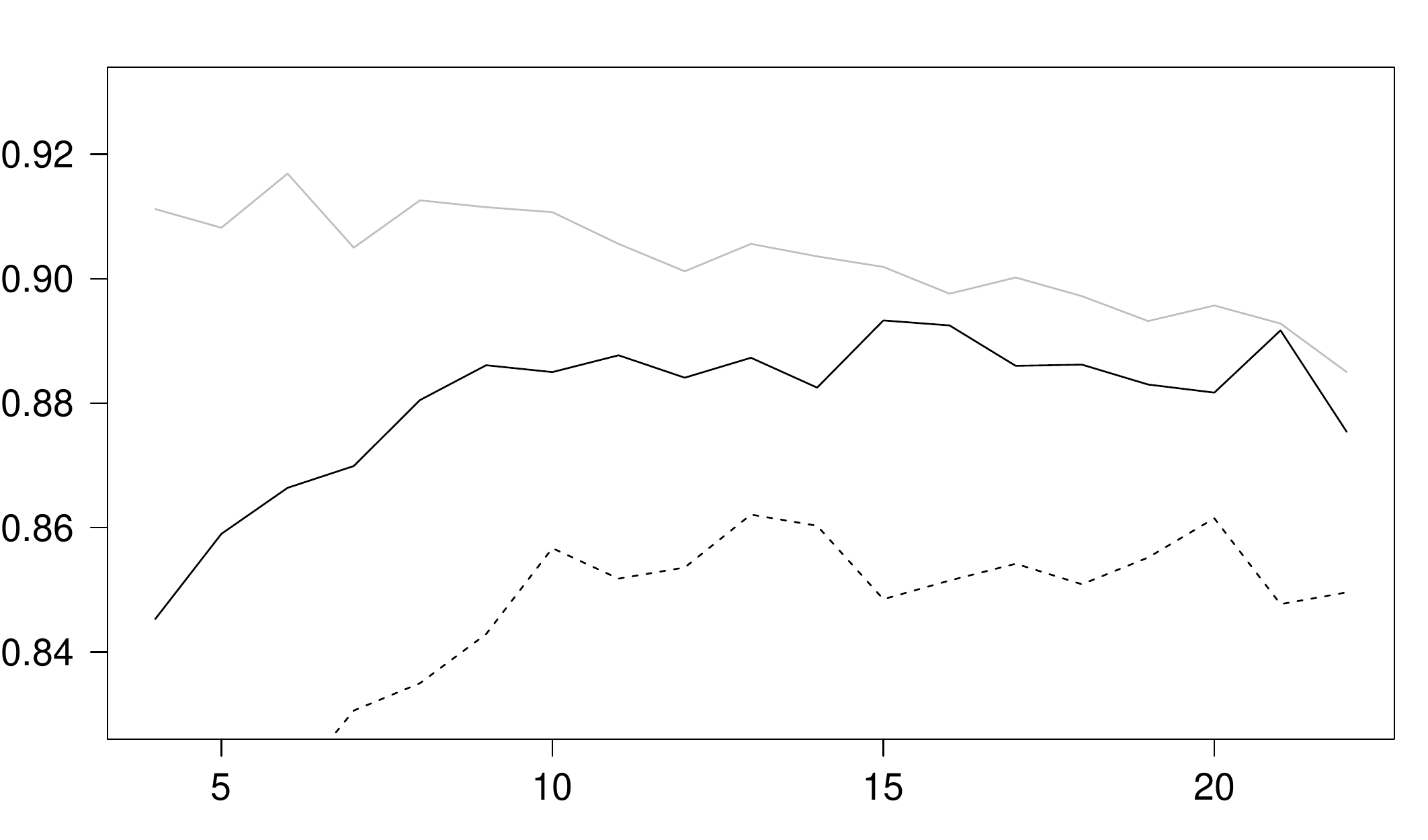}
        \caption{Simulated coverage probabilities for plug-in bootstrap (grey), new bootstrap (black), subsampling (dashed) as a function of the block length $l$, sample size $n=200$, AR-coefficient $\alpha=0.4$}
\end{figure}

\section{Proofs of the theorems}
\subsection{Preliminary results}
In this section we formulate some necessary results that will be used in the proofs of our theorems. 
\begin{lemma}\label{lem1}(Dehling,Wendler \cite{deh2})
Let $\left(X_n\right)_{n\in\mathds{N}}$ be a stationary sequence and $h(x,y)$ a degenerate kernel satisfying $\iint\left|h\left(x_{1},x_{2}\right)\right|^{12+\delta}dF\left(x_{1}\right)dF\left(x_{2}\right)\leq M$ for some $\delta>0$ and $E\left|h\left(X_{1},X_{1+m}\right)\right|^{4+\frac{2}{3}\delta} \leq M$ for all $m\in\N$. Let be $\tau\geq0$ such that one of the following two conditions holds:
\begin{enumerate}
 \item  $\left(X_n\right)_{n\in\mathds{N}}$ is absolutely regular and $\sum_{k=0}^{n}k\beta^{\frac{\delta}{2+\delta}}\left(k\right)=O\left(n^\tau\right)$.
\item $\left(X_n\right)_{n\in\mathds{N}}$ is strongly mixing, $E\left|X_1\right|^\gamma<\infty$ for a $\gamma>0$, $h$ is $\mathcal{P}$-Lipschitz-continuous with constant $C>0$ and $\sum_{k=0}^{n}k\alpha^{\frac{2\gamma\delta}{\gamma\delta+\delta+5\gamma+2}}\left(k\right)=O\left(n^\tau\right)$.
\end{enumerate}
Then:
\begin{equation*}
\sum_{i_{1},i_{2},i_{3},i_{4}=1}^{n}\left|E\left[h_{2}\left(X_{i_{1}},X_{i_{2}}\right)h_{2}\left(X_{i_{3}},X_{i_{4}}\right)\right]\right|=O\left(n^{2+\tau}\right),
\end{equation*}
\end{lemma}
\begin{lemma}\label{lem2}(Yoshihara \cite{yosh})
  Let $\left(X_n\right)_{n\in\mathds{N}}$ be a stationary sequence of absolutely regular random variables. Assume that the  the following conditions hold (for some $\delta >0$) and
\begin{equation*}
\iint\left|h\left(x_{1},x_{2}\right)\right|^{4+\delta}dF\left(x_{1}\right)dF\left(x_{2}\right)\leq M,
\end{equation*}
\begin{equation*}
E\left|h\left(X_{1},X_{1+k}\right)\right|^{4+\delta} \leq M \text{ for all } k\in\N
\end{equation*}
\begin{equation*}
\beta\left(n\right)=O\left(n^{-\frac{3(4+\delta')}{\delta'}}\right)\ \text{ for some } 0<\delta'< \delta.
\end{equation*}
Then
\begin{equation*}
 E\left(n^{2}U^{4}_{n}(h_{2})\right)=O(n^{-1-\eta}),\ \ \eta=\min \left\{ 12\dfrac{\delta-\delta'}{\delta'(4+\delta)},1\right\} 
\end{equation*}
\end{lemma}
\begin{lemma}\label{lem3} (Yokoyama \cite{yoko}) Let $\left(X_n\right)_{n\in\N}$ be a stationary strongly mixing sequence of random variables with $E X_1 = \mu$ and $\left\|X_1\right\|_{2 + \delta} < \infty$ for some $ 0 < \delta \leq \infty$. Suppose that $2\leq s < 2 + \delta$ and
\[
\sum^\infty_{n =1} n^{\frac{s}{2} - 1} \left(\alpha (n)\right)^{(2 + \delta - s)/ (2 + \delta)} < \infty.
\]
Then there exists a constant $C$ depending only on $s, \delta$ and the mixing coefficients $\left(\alpha(n)\right)_{n\in\N}$ such that
\begin{equation*}
E \left| \sum^n_{i =1} (X_i - \mu)\right|^s \leq C n^{s/2}.
\end{equation*}
\end{lemma}

\subsection{Proof of Theorem \ref{theo1}}
Let us introduce blocks of indices
\begin{align*}
B_1(i)&= \left\{i, \ldots, i+l-1\right\},\quad i = 1, \ldots, n , \\
B_2(i)&=\left\{(i-1)l+1, \ldots, il\right\},\quad i = 1, \ldots, m.
\end{align*}
Note that the blocks $\lbrace B_{1}(i),i=1,...,n\rbrace$ and $\lbrace B_{2}(i),i=1,...,m \rbrace$ correspond to the circular and nonoverlapping blocking methods, respectively. In the circular blocking method instead of the sample $ X_{1},...,X_{n}$ we consider the completed sample $ X_{1},...,X_{n},X_{1},...,X_{l-1}.$ A fixed realization of the sample $ X_{1},...,X_{n}$ we denote by $ x_{1},...,x_{n}$ and set
\begin{align*}
\overline{x}_{1,1}(i)=\frac{1}{l}\sum_{k\in B_{1}(i)} x_{k}  , i=1,...,n , \\ 
\overline{x}_{2,1}(i)=\frac{1}{l}\sum_{k\in B_{2}(i)} x_{k}   , i=1,...,m , \displaybreak[0]\\
X_{j}^{1*}=\sum_{i=1}^{n}I(T_{1}(j)=i)\overline{x}_{1,1}(i) ,\\
X_{j}^{2*}=\sum_{i=1}^{m}I(T_{2}(j)=i)\overline{x}_{2,1}(i) ,\displaybreak[0]\\
\overline{X}_{n,m}^{1*}=\frac{1}{m}\sum_{i=1}^{m}X_{i}^{1*} , \\
\overline{X}_{n,m}^{2*}=\frac{1}{m}\sum_{i=1}^{m}X_{i}^{2*}.
\end{align*}

 First we will prove the statement of the theorem for $i=1$ (circular bootstrap). Using Hoeffding decomposition we have
\begin{align*}
&\left(U^{1\star}_{n}\left(h\right)-E^{\star}\left[U^{1\star}_{n}\right]\right)=\frac{1}{m} \sum^m_{j =1}\left(u^{1*}_{j}-E^{\star}U^{1*}_{n}(h)\right)\\
= \frac{1}{m}& \sum^m_{j=1}\bigg[\theta+\frac{2}{l}\sum_{k \in B_{1}(T_{1}(j))}h_{1}(x_{k})+\frac{2}{l(l-1)}\sum_{(a,b) \in B^{'}_{1}(T_{1}(j))}h_{2}(x_{a},x_{b})\\
&-\big(\theta+\frac{1}{n}\sum^n_{i=1}\frac{2}{l}\sum_{k \in B_{1}(i)}h_{1}(x_{k})+\frac{1}{n} \sum^n_{i=1}\frac{2}{l(l-1)}\sum_{(a,b) \in B^{'}_{1}(i)}h_{2}(x_{a},x_{b})\big)\bigg] \displaybreak[0]\\
=\frac{1}{m}& \sum^m_{j=1}\bigg[\frac{2}{l}\big(\sum_{k \in B_{1}(T_{1}(j))}h_{1}(x_{k})-\frac{1}{n}\sum^n_{r =1}h_{1}(x_{r})\big)\\
&+\frac{2}{l(l-1)}\big(\sum_{(a,b) \in B^{'}_{1}(T_{1}(j))}h_{2}(x_{a},x_{b})-\frac{1}{n} \sum^n_{r=1}\frac{2}{l(l-1)}\sum_{(a,b) \in B^{'}_{1}(r)}h_{2}(x_{a},x_{b})\big)\bigg]\displaybreak[0]\\
= \frac{2}{ml}& \sum^m_{j =1}\sum_{k \in B_{1}(T_{1}(j))}\big( h_{1}(x_{k})-\frac{1}{n}\sum^n_{r =1}h_{1}(x_{r})\big)\\
&+\frac{2}{ml(l-1)}\sum^m_{j =1}\sum_{(a,b) \in B^{'}_{1}(T_{1}(j))}\big(h_{2}(x_{a},x_{b})-\frac{1}{n} \sum^n_{r=1}\sum_{(a,b) \in B^{'}_{1}(r)}h_{2}(x_{a},x_{b})\big)
\end{align*}
Using the latter
\begin{align*}
&\sqrt{ml}\left(U^{1\star}_{n}\left(h\right)-E^{\star}\left[U^{1\star}_{n}\right]\right)=\frac{2}{\sqrt{ml}} \sum^m_{j =1}\Big[\sum_{k \in B_{1}(T_{1}(j))}h_{1}(x_{k})-\frac{1}{n}\sum^n_{i =1}h_{1}(x_{i})\Big]\\
+&\frac{2}{\sqrt{ml}(l-1)}\sum^m_{j =1}\Big[\sum_{(a,b) \in B^{'}_{1}(T_{1}(j))}(h_{2}(x_{a},x_{b})-\frac{1}{n} \sum^n_{r=1}\sum_{(a,b) \in B^{'}_{1}(r)}h_{2}(x_{a},x_{b}))\Big]\\
=&I_{n}+I\! I_{n}.
\end{align*}
In order to prove the first statement of Theorem \ref{theo1}, we first show that $I\! I_{n}\rightarrow 0 $ as $n\rightarrow\infty$. Note that
\begin{align*}
\Var^{\star}I\! I_{n}=\frac{4}{nl(l\!-\!1)^{2}}\sum^n_{i =1}\bigg(\sum_{(a,b) \in B^{'}_{1}(i)}\Big(h_{2}(x_{a},x_{b})-\frac{1}{n} \sum^n_{i =1}\sum_{(a,b) \in B^{'}_{1}(i)}h_{2}(x_{a},x_{b})\Big)\bigg)^{2}\\
=\frac{1}{n}\sum^n_{i =1}\Big(\frac{2}{\sqrt{l}(l\!-\!1)}\sum_{(a,b) \in B^{'}_{1}(i)}h_{2}(x_{a},x_{b})\Big)^{2}-\Big(\frac{1}{n} \sum^n_{i =1}\frac{2}{\sqrt{l}(l\!-\!1)}\sum_{(a,b) \in B^{'}_{1}(i)}h_{2}(x_{a},x_{b})\Big)^{2}.
\end{align*}
By Lemma \ref{lem1} we have as $n\rightarrow\infty$
\begin{equation*}
 E\Big(\frac{2}{\sqrt{l}(l-1)}\sum_{(a,b) \in B^{'}_{1}(i)}h_{2}(X_{a},X_{b})\Big)^2\rightarrow 0
\end{equation*}
and consequently $\Var^{\star}I\! I_{n}\rightarrow 0$ in probability as $n\rightarrow\infty$. With the Chebyshev inequality, it follows that $I\! I_{n}\rightarrow 0 $ in conditional probability.

It remains to show the convergence of $I_n$. By Theorem 3.2 of Lahiri \cite{lahi} we have as $n\rightarrow\infty$
\begin{equation*}
\sup_{x\in\R} \bigg| P^{\star}\big(I_{n} \leq x \big) -P \Big( \frac{2}{\sqrt{n}}\sum^n_{i =1}h(X_{i}) \leq x\Big)\bigg| \rightarrow 0.
\end{equation*}
in probability. Furthermore, we have $\Var\left(\sqrt{n}U_{n}(h_{2})\right)\rightarrow 0$ by Lemma \ref{lem1} and can conclude that
\begin{equation*}
\sup_{x\in\R} \bigg|P \Big( \frac{2}{\sqrt{n}}\sum^n_{i =1}h(X_{i}) \leq x\Big)-P\Big(\sqrt{n}(U_n(h)-\theta)\leq x\Big)\bigg| \rightarrow 0.
\end{equation*}
Using Slutzky's lemma and  $I\! I_{n}\rightarrow 0 $ in probability, we arrive at at the first statement of the Theorem:
\begin{equation*}
 \sup_{x\in\mathds{R}}\left|P^{\star}\left[\sqrt{ml}\left(U^{i\star}_{n}\left(h\right)-E_{\star}\left[U^{i\star}_{n}\right]\right)\leq x\right]-P\left[\sqrt{n}\left(U_{n}\left(h\right)-\theta\right)\leq x\right]\right|\rightarrow0.
\end{equation*}
For the second statement, we make use of Theorem 3.2 of Lahiri \cite{lahi} again, which also states that
\begin{equation*}
\left|\Var^{\star}I_{n}-\Var\Big(\frac{2}{\sqrt{n}}\sum^n_{i =1}h(X_{i})\Big)\right|\rightarrow 0 
\end{equation*}
 This together with $\Var^{\star}I\! I_{n}\rightarrow 0$ in probability and $\Var\left(\sqrt{n}U_{n}(h_{2})\right)\rightarrow 0$ (see Lemma \ref{lem1}) leads to
\begin{equation*}
\left|\Var^\star\left[\sqrt{ml}U_n^{i\star}\left(h\right)\right]-\Var\left[\sqrt{n}U_n\left(h\right)\right]\right|\rightarrow0
\end{equation*}
in probability as $n\rightarrow\infty$.

In the case of $i=2$ (nonoverlapping bootstrap) analogously to previous one we have
\begin{multline*}
\sqrt{ml}\left(U^{2\star}_{n}\left(h\right)-E_{\star}\left[U^{2\star}_{n}\right]\right)\\
=\frac{2}{\sqrt{ml}} \sum^m_{j =1}\sum^m_{i =1}I(T_{2}(j)=i)\sum_{k \in B_{1}(i)}\big(h_{1}(x_{k})-\frac{1}{ml}\sum^{ml}_{i =1}h_{1}(x_{i})\big)\\
+\frac{2}{\sqrt{ml}(l-1)}\sum^m_{j =1}\sum_{(a,b) \in B^{'}_{1}(T_{2}(j))}\big(h_{2}(x_{a},x_{b})-\frac{1}{m} \sum^m_{i =1}\sum_{(a,b) \in B^{'}_{1}(i)}h_{2}(x_{a},x_{b})\big)=\\
=I^{'}_{n}+I\! I^{'}_{n}.
\end{multline*}
Theorem 3.2 of Lahiri \cite{lahi} implies as $n\rightarrow\infty$
\begin{equation*}
\sup_{x\in\R} \left| P^{\star}\left(I^{'}_{n} \leq x \right) -P \left( \frac{2}{\sqrt{n}}\sum^n_{i =1}h(X_{i}) \leq x\right)\right| \rightarrow 0
\end{equation*}
\begin{equation*}
\left|\Var^{\star}I^{'}_{n}-\Var\left(\frac{2}{\sqrt{n}}\sum^n_{i =1}h(X_{i})\right)\right|\rightarrow 0
\end{equation*}
in probability as $n\rightarrow\infty$. It remains to show that $\Var^{\star}I\! I^{'}_{n}\rightarrow 0$ in probability, the rest of the proof can be done along the lines of the case $i=1$ (circular bootstrap). We have
\begin{align*}
\Var^{\star}I\! I^{'}_{n}=\frac{4}{ml(l\!-\!1)^{2}}\sum^m_{i =1}\left(\sum_{(a,b) \in B^{'}_{2}(i)}\big(h_{2}(x_{a},x_{b})-\frac{1}{m} \sum^m_{i =1}\sum_{(a,b) \in B^{'}_{2}(i)}h_{2}(x_{a},x_{b})\big)\right)^{2}\\
=\frac{1}{m}\sum^m_{i =1}\big(\frac{2}{\sqrt{l}(l\!-\!1)}\sum_{(a,b) \in B^{'}_{2}(i)}h_{2}(x_{a},x_{b})\big)^{2}-\big(\frac{1}{m} \sum^m_{i =1}\frac{2}{\sqrt{l}(l\!-\!1)}\sum_{(a,b) \in B^{'}_{2}(i)}h_{2}(x_{a},x_{b})\big)^{2}.
\end{align*}
By Lemma \ref{lem1} we have as $n\rightarrow\infty$
\begin{equation*}
E\left(\frac{2}{\sqrt{l}(l-1)}\sum_{(a,b) \in B^{'}_{2}(i)}h_{2}(X_{a},X_{b})\right)^2\rightarrow 0
\end{equation*}
so $E\left(\Var^{\star}I\! I^{'}_{n}\right)\rightarrow 0$.

\subsection{Proof of Theorem \ref{theo2}}
Theorem \ref{theo2} can be proved in the same way as Theorem \ref{theo1} and therefore the proof is omitted.

\subsection{Proof of Theorem \ref{theo3}} To shorten the proof, we concentrate on the case $k=1$ (circular bootstrap). We split the mean squared error into three parts
\begin{multline*}
E\left(\Var^{\star}\big(\sqrt{ml}U_{n}^{k*}(h)\big)-\Var\big(\sqrt{n}U_{n}(h)\big)\right)^{2}\\
\leq 3 E\Big(\Var^{\star}\big(\frac{1}{\sqrt{ml}}\sum_{i=1}^{ml}h_1(X_{1,i}^\star)\big)-\Var\big(\frac{1}{\sqrt{n}}\sum_{i=1}^nh_1(X_i)\big)\Big)^{2}\\
+3 E\Big(\Var^{\star}\big(\frac{1}{\sqrt{ml}}\sum_{i=1}^{ml}h_1(X_{1,i}^\star)\big)-\Var^{\star}\big(\sqrt{ml}U_{n}^{k*}(h)\big)\Big)^{2}\\
+3 E\Big(\Var(\sqrt{n}U_{n}(h)))-\Var(\frac{1}{\sqrt{n}}\sum_{i=1}^nh_1(X_i))\Big)^{2}=3I_n+3I\! I_n+3I\! I\! I_n
\end{multline*}
For the first summand, we know by Theorem 3.1 of Politis and White \cite{pol2} that
\begin{equation*}
I_n= E\Big(\Var^{\star}\big(\frac{1}{\sqrt{ml}}\sum_{i=1}^{ml}h_1(X_{1,i}^\star)\big)-\Var\big(\frac{1}{\sqrt{n}}\sum_{i=1}^nh_1(X_i)\big)\Big)^{2}=O(\frac{1}{l^2}+\frac{l}{n}).
\end{equation*}
For the second summand, we make use of the equation $a^2-b^2=(a+b)(a-b)$ and the H\"older-inequality (first for the bootstrap expectation and then for the unconditional expectation) to obtain
\begin{multline*}
I\! I_n=E\Big(\Var^{\star}(\frac{1}{\sqrt{ml}}\sum_{i=1}^{ml}h_1(X_{1,i}^\star))-\Var^{\star}(\sqrt{ml}U_{n}^{k*}(h))\Big)^{2}\\
=E\Bigg(E^{\star}\bigg(\Big(\frac{1}{\sqrt{ml}}\sum_{i=1}^{ml}h_1(X_{1,i}^\star)-E^\star\frac{1}{\sqrt{ml}}\sum_{i=1}^{ml}h_1(X_{1,i}^\star)-\sqrt{ml}U_{n}^{k*}(h)+E^\star\sqrt{ml}U_{n}^{k*}(h)\Big)\\
\shoveright{\Big(\frac{1}{\sqrt{ml}}\sum_{i=1}^{ml}h_1(X_{1,i}^\star)-E^\star\frac{1}{\sqrt{ml}}\sum_{i=1}^{ml}h_1(X_{1,i}^\star)+\sqrt{ml}U_{n}^{k*}(h)-E^\star\sqrt{ml}U_{n}^{k*}(h)\Big)\bigg)\Bigg)^{2}}\displaybreak[0]\\
\leq E\bigg(E^{\star}\Big(\frac{1}{\sqrt{ml}}\sum_{i=1}^{ml}h_1(X_{1,i}^\star)-E^\star\frac{1}{\sqrt{ml}}\sum_{i=1}^{ml}h_1(X_{1,i}^\star)-\sqrt{ml}U_{n}^{k*}(h)+E^\star\sqrt{ml}U_{n}^{k*}(h)\Big)^2\\
\shoveright{E^{\star}\Big(\frac{1}{\sqrt{ml}}\sum_{i=1}^{ml}h_1(X_{1,i}^\star)-E^\star\frac{1}{\sqrt{ml}}\sum_{i=1}^{ml}h_1(X_{1,i}^\star)+\sqrt{ml}U_{n}^{k*}(h)-E^\star\sqrt{ml}U_{n}^{k*}(h)\Big)^{2}\bigg)}\displaybreak[0]\\
\shoveleft{= E\bigg(\Var^{\star}\Big(\frac{1}{\sqrt{ml}(l-1)}\sum_{j=1}^m\sum_{i=1}^nI(T_1(j)=i)\sum_{(a,b)\in B_1'(i)}h_2(X_a,X_b)\Big)}\\
\Var^{\star}\Big(2\frac{1}{\sqrt{ml}}\sum_{i=1}^{ml}h_1(X_{1,i}^\star)+\frac{1}{\sqrt{ml}(l-1)}\sum_{j=1}^m\sum_{i=1}^nI(T_1(j)=i)\sum_{(a,b)\in B_1'(i)}h_2(X_a,X_b)\Big)\bigg)\displaybreak[0]\\
\leq \bigg(E\Big(\Var^{\star}\big(\frac{1}{\sqrt{ml}(l-1)}\sum_{j=1}^m\sum_{i=1}^nI(T_1(j)=i)\sum_{(a,b)\in B_1'(i)}h_2(X_a,X_b)\big)\Big)^2\bigg)^{\frac{1}{2}}\\
\bigg(E\Big(\Var^{\star}\big(2\frac{1}{\sqrt{ml}}\sum_{i=1}^{ml}h_1(X_{1,i}^\star)+\frac{1}{\sqrt{ml}(l-1)}\sum_{j=1}^m\sum_{(a,b)\in B_1'(T_1(j))}h_2(X_a,X_b)\big)\Big)^2\bigg)^{\frac{1}{2}}\\
=A_n\times B_n.
\end{multline*}
We will treat these two factors separately, starting with $A_n$. By the definition of the bootstrap procedure and the stationarity of the sequence, we get the following:
\begin{multline*}
A_n=\bigg(E\Big(\Var^{\star}\big(\frac{1}{\sqrt{ml}(l-1)}\sum_{j=1}^m\sum_{i=1}^nI(T_1(j)=i)\sum_{(a,b)\in B_1'(i)}h_2(X_a,X_b)\big)\Big)^2\bigg)^{\frac{1}{2}}\\
\leq \bigg(E E^{\star}\Big(\frac{1}{\sqrt{ml}(l-1)}\sum_{j=1}^m\sum_{i=1}^nI(T_1(j)=i)\sum_{(a,b)\in B_1'(i)}h_2(X_a,X_b)\Big)^4\bigg)^{\frac{1}{2}}\\
=\bigg(E \frac{1}{l^2(l-1)^4}\frac{1}{n}\sum_{i=1}^n\Big(\sum_{(a,b)\in B_1'(i)}h_2(X_a,X_b)\Big)^4\bigg)^{\frac{1}{2}}\\
=\bigg( \frac{1}{l^2(l-1)^4}E\Big(\sum_{(a,b)\in B_1'(i)}h_2(X_a,X_b)\Big)^4\bigg)^{\frac{1}{2}}=O(\frac{1}{l}),
\end{multline*}
as we know from Lemma 3 of Yoshihara \cite{yosh} that
\begin{equation*}
E\Big(\sum_{1\leq i<j\leq k}h_2(X_i,X_j)\Big)^4=O(l^4).
\end{equation*}
With the help of the inequality $(a+b)^2\leq 2 a^2+2b^2$, we split the second summand into two parts:
\begin{multline*}
B_n\\
=\bigg(E\Big(\Var^{\star}\big(2\frac{1}{\sqrt{ml}}\sum_{i=1}^{ml}h_1(X_{1,i}^\star)+\frac{1}{\sqrt{ml}(l-1)}\sum_{j=1}^m\sum_{(a,b)\in B_1'(T_1(j))}h_2(X_a,X_b)\big)\Big)^2\bigg)^{\frac{1}{2}}\\
\leq 16 \bigg(E\Big(\Var^{\star}\big(\frac{1}{\sqrt{ml}}\sum_{i=1}^{ml}h_1(X_{1,i}^\star)\big)\Big)^2\bigg)^{\frac{1}{2}}+4A_n=16B'_n+4A_n,
\end{multline*}
so it suffices to study $B'_n$. We use again stationarity and the definition of the bootstrap:
\begin{multline*}
B'_n=\bigg(E\Big(\Var^{\star}\big(\frac{1}{\sqrt{ml}}\sum_{i=1}^{ml}h_1(X_{1,i}^\star)\big)\Big)^2\bigg)^{\frac{1}{2}}=\bigg(E\Big(\Var^{\star}\big(\frac{1}{\sqrt{l}}\sum_{i=1}^{l}h_1(X_{1,i}^\star)\big)\Big)^2\bigg)^{\frac{1}{2}}\displaybreak[0]\\
\leq \bigg(EE^{\star}\big(\frac{1}{\sqrt{l}}\sum_{i=1}^{l}h_1(X_{1,i}^\star)\big)^4\bigg)^{\frac{1}{2}}=\bigg(E\frac{1}{n}\sum_{j=1}^n\big(\frac{1}{\sqrt{l}}\sum_{i=j}^{j+l-1}h_1(X_{i})\big)^4\bigg)^{\frac{1}{2}}\\
=\bigg(E\big(\frac{1}{\sqrt{l}}\sum_{i=j}^{j+l-1}h_1(X_{i})\big)^4\bigg)^{\frac{1}{2}}\leq C.
\end{multline*}
The last inequality follows from Lemma \ref{lem3}. This shows that $I\! I_n=O(\frac{1}{l})$. In the same way one can show that $I\! I\! I_n$ is of the same order, which completes the proof.

\end{document}